    \input amstex
     \documentstyle{amsppt}
   \pageno=1
\NoRunningHeads

     \magnification = 1200
    \pagewidth{5.3 true in}
    \pageheight{8.5 true in}
   \hoffset .6 true in

   \NoBlackBoxes
     \parskip=3 pt
     \parindent = 0.3 true in

\def\({\bigl(}
\def\){\bigr)}
\def\Z{\Bbb Z}

\def\Q{\Bbb Q}
\def\R{\Bbb R}
\def\C{\Bbb C}
\def\FF{\Bbb F}
\def\F{\Cal F}

\def\G{\Cal G}
\def\R{\Cal R}
\def\S{\Cal S}

\def\N{\operatorname{N}}
\def\Tr{\operatornamewithlimits{Tr}}
\def\Oc{\Cal O}

\def\fa{\frak a}

\def\p{\frak p}

\def\q{\frak q}

\def\Fit{{\operatornamewithlimits{Fit}}}
\def\Ann{{\operatornamewithlimits{Ann}}}

\def\Gal{{\operatorname{Gal}}}

\def\coker{{\operatorname{coker}}}

\def\Stick{{\operatornamewithlimits{Stick}}}
\def\op{{{\bigoplus}}}

\def\FI{{\Fit_{\E/F}^S}(1)}
\def\SI{{\Stick_{\E/F}^S}(-1)}
\def\FIb{{\Fit_{E/F}^S}(1)}
\def\SIb{{\Stick_{E/F}^S}(-1)}

\def\Rb{{\overline{\R}}}
\def\Gb{{\overline{G}}}
\def\gb{{\overline{\gamma}}}
\def\Ghat{{\hat{G}}}
\def\E{{\Cal E}}
\def\S{{\Cal S}}

\def\OE{\Oc_E^S}
\def\OEc{\Oc_\E^S}
\def\OEchi{\Oc_{E_\chi}^S}
\def\OF{\Oc_F^S}
\def\KE{K_2(\OE)}
\def\KEchi{K_2(\Oc_{E_\chi}^S)}
\def\kechi{k_2^S(E_\chi)}
\def\KEc{K_2(\OEc)}
\def\KF{K_2(\OF)}
\def\kf{k_2^S(F)}
\def\WE{W_2(E)}
\def\wechi{w_2(E_\chi)}
\def\WEc{W_2(\E)}
\def\WF{W_2(F)}
\def\thE{\theta_{\E/F}^S(-1)}
\def\Ba{1}
\def\Co{2}
\def\JTT{20}
\def\Ga{4}
\def\QuFG{13}
\def\MW{11}
\def\Wi{22}
\def\Ko{8}
\def\Kol{7}
\def\KoMo{9}
\def\Ka{6}
\def\SaSi{16}
\def\Sa{14}
\def\San{15}
\def\So{19}
\def\Gr{5}
\def\KS{{17}}
\def\DR{{3}}
\def\Sn{{18}}
\def\Ta{{21}}
\def\Mi{{12}}
\def\Les{{10}}
    \topmatter
    \title
Values at $s=-1$ of $L$-functions \\
for multi-quadratic extensions \\
of number fields, and \\
the Fitting ideal of the tame kernel
    \endtitle
\author
Jonathan W. Sands
\\
Dept. of Mathematics and Statistics
\\
University of Vermont, Burlington VT 05401 USA
\\
Jonathan.Sands\@uvm.edu
\endauthor

\thanks
2000 {\it MSC:} Primary 11R42; Secondary 11R70, 19F27.
\endthanks

\abstract
  Fix a Galois extension $\E/F$ of totally real number fields such that the
Galois group $G$ has exponent 2. Let $S$ be a finite set of primes
of $F$ containing the infinite primes and all those which ramify in
$\E$, let $S_\E$ denote the primes of $\E$ lying above those in $S$,
and let $\Oc_\E^S$ denote the ring of $S_\E$-integers of $\E$. We
then compare the Fitting ideal of $K_2(\Oc_{\E}^S)$ as a
$\Z[G]$-module with a higher Stickelberger ideal.
The two extend to the same ideal in the maximal order of $\Q[G]$,
and hence in $\Z[1/2][G]$. Results in $\Z[G]$
are obtained under the assumption of the Birch-Tate
conjecture, especially for biquadratic extensions,
where we compute the index of the higher Stickelberger ideal.
 We find a sufficient
condition for the Fitting ideal to contain the higher
Stickelberger ideal in the case where $\E$ is a biquadratic
extension of $F$ containing the first layer of the cyclotomic
$\Z_2$-extension of $F$, and describe a class of biquadratic
extensions of $F=\Q$ that satisfy this condition.
\endabstract

    \endtopmatter

     \document
\baselineskip=20pt

\bigskip\bigskip\bigskip
 \vfill\eject

\heading{I. Introduction}
\endheading

 Fix an abelian Galois
extension of number fields $\E/F$ and let $G$ be the Galois group. Also
fix a finite set $S$ of primes of $F$ which contains all of the
infinite primes of $F$ and all of the primes which ramify in $\E$.
Associated with this data is a Stickelberger function, or
equivariant $L$-function, $\theta_{\E/F}^S(s)$. It is a meromorphic
function of $s$ with values in the group ring $\C[G]$. To define it,
let $\p$ run through the (finite) primes of $F$ not in $S$, and
$\fa$ run through integral ideals of $F$ which are relatively prime
to each of the elements of $S$. Also let ${\N} \fa$ denote the absolute
norm of the ideal $\fa$ and
 $\sigma_{\fa} \in G$ denote the well-defined automorphism
attached to $\fa$ via the Artin map. Then

$${
\theta_{\E/F}^S(s)=
\sum\Sb {\frak a} \ \text{integral} \\ ({\frak a},S)=1 \endSb
\dfrac{1}{ {\N} {\frak a}^s}\sigma_{\fa}^{-1}=
\prod\Sb \text{prime } \frak p \notin S \endSb
\bigl(1-\dfrac{1}{ {\N} {\frak p}^{s}}} \sigma_{\frak p}^{-1}\bigr)^{-1}.$$

 These expressions converge for the real part of $s$ greater than 1, and the function they define
extends meromorphically to all of $\C$. When $\E=F$, the function $\theta_{F/F}^S(s)$ is simply the
identity automorphism of $F$ times $\zeta_F^S(s)$, the Dedekind zeta-function of $F$ with Euler factors
for the primes in $S$ removed.

The function $\theta_{\E/F}^S(s)$ is connected with the arithmetic of
the number fields $\E$ and $F$ in ways one would like to make as
precise as possible.
The ring of $S$-integers $\OF$ of $F$ is defined to
be the set of elements of $F$ whose valuation is non-negative at
every prime not in $S$. Similarly, define the ring $\OEc$ of
$S$-integers of $\E$ to be the set of elements of $E$ whose valuation
is non-negative at every prime not in $S_\E$, the set of all primes
of $\E$ which lie above some prime in $S$. The function $\zeta_F^S(s)$
may be viewed as the zeta-function of the Dedekind domain $\OF$.

We are interested in the ``higher Stickelberger element" $\theta_{\E/F}^S(-1)$,
which lies in $\Q[G]$ by the theorem of Klingen-Siegel [\KS],
and is related to the algebraic $K$-group $\KEc$. This
group is known to be finite by [\Ga] and [\QuFG], and could be called the
$S$-tame kernel of $\E$. It contains the tame kernel $K_2(\Oc_\E)$ as a subgroup.
Another piece of the arithmetic interpretation
of $\theta_{\E/F}^S(-1)$ involves a group of roots of unity.
Let $\mu_\infty$ denote the group of all roots of unity in an algebraic closure
$\overline{\Q}$ of $\Q$ containing $\E$, and let $\G$ denote the
Galois group of $\overline{\Q}/\Q$. Define $W_2=W_2(\overline{\Q})$
to be the $\Z[\G]$-module whose underlying group is $\mu_\infty$,
with the action of $\gamma \in \G$ on $\omega \in W_2$ given by
$\omega^\gamma=\gamma^2(\omega)$. For any subfield $L$ of
$\overline{\Q}$, let $W_2(L)$ be the submodule fixed under this
action by the Galois group of $\overline{\Q}$ over $L$. Then
$W_2(\E)$ naturally becomes a $\Z[G]$-module, where the action of $G$
arises by lifting elements of $G$ to $\G$ and then using the action
of $\G$ just defined. One easily sees that the $G$-fixed
submodule $W_2(\E)^G$ equals $W_2(F)$. We use the notation
$w_2(L)=|W_2(L)|$, which we note is finite for any algebraic
number field $L$.

Vital to our approach is the conjecture of Birch and Tate (see
section 4 of [\JTT]), which gives a precise arithmetic
interpretation of $\zeta_F^S(-1)$. We state a form of it
for an arbitrary finite set $S$ which is easily seen to be equivalent to the original
conjecture for the minimal choice of the set $S$, containing just the infinite primes
 (see Corollary 3.3 of [\SaSi]).

\proclaim{Conjecture 1.1 (Birch-Tate)}
Suppose that $F$ is totally real and the finite set $S$ contains
the infinite primes of $F$. Then
$$\zeta_F^S(-1)= (-1)^{|S|} \dfrac{|\KF|}{w_2(F)}$$
\endproclaim

 Deep results on Iwasawa's Main conjecture in [\MW] and [\Wi] lead to the following
(see [\Ko]).

\proclaim{Proposition 1.2}
The Birch-Tate Conjecture holds if $F$ is abelian over $\Q$, and the odd part
holds for all totally real $F$.
\endproclaim

Kolster [\Kol] has shown that the 2-part of the Birch-Tate conjecture for $F$ would
follow from the 2-part of Iwasawa's Main conjecture for $F$.

For any module $M$ over a ring $R$, we let $\Ann_R(M)$ denote the annihilator of
$M$ in $R$. The following result is proved in [\San, Thm 1.3].

\proclaim{Proposition 1.3} Let $E/F$ be a relative quadratic extension of
totally real number fields, with Galois group $\Gb$. Let $S$ contain the infinite primes and
those which ramify in $E/F$.
Assume that the 2-part of the Birch-Tate conjecture holds for $E$ and for $F$.
Then the (first) Fitting ideal of $\KE$ as a $\Z[\Gb]$-module is
$$\Fit_{\Z[\Gb]}(\KE)=\Ann_{\Z[\Gb]}(W_2(E))\theta_{E/F}^S(-1).$$
More specifically, this ideal equals its extension to the maximal
order of $\Q[\Gb]$ if and only if it is not principal, and this
happens exactly when $E$ is not the first layer of the cyclotomic
$\Z_2$-extension of $F$. Without the assumption of the Birch-Tate
conjecture, the ideals $\Fit_{\Z[\Gb]}(\KE)$ and
$\Ann_{\Z[\Gb]}(W_2(E))\theta_{E/F}^S(-1)$ have the same extension
to $\Z[1/2][\Gb]$.
\endproclaim

In this paper, we build upon Proposition 1.3, obtaining more general
results which suggest a close relationship
between $\FI=\Fit_{\Z[G]}(\KEc)$ and the higher
Stickelberger ideal $\SI=\Ann_{\Z[G]}(W_2(\E))\theta_{\E/F}^S(-1)$, particularly
when $G$ has exponent 2. Here the theorem of Deligne and Ribet [\DR] guarantees that $\SI$
is an ideal in the integral group ring $\Z[G]$.

Part of the motivation for this work is to compare the situation concerning
$\KEc$ and $\thE$, with that of the $S$-ideal class group of $\E$ and
the Stickelberger element $\theta_{\E/F}^S(0)$. Non-triviality in the latter situation
requires that $\E$ be a CM field and $F$ be totally real. Then under the assumption of the Equivariant Tamagawa Number Conjecture,
 Greither [\Gr] has recently obtained an equality of $p$-parts for odd primes $p$ between an ideal built from Stickelberger
 elements, and the Fitting ideal of the Pontrjagin dual
 of the minus part of the ideal class group of $\E$. Thus one might wonder whether
 the Pontrjagin dual should occur in our investigation. So far, it does not.
 Taking the dual does not affect
 the Fitting ideal in the case of a cyclic group $G$, so the case of a non-cyclic group
 is of special interest to us. Our results tend to demonstrate the naturality of
 not taking the dual. For instance, the property
 of the Fitting ideal of the $S$-tame kernel with respect to change of extension field
 in Proposition 2.4 seems to have no straightforward counterpart if the $S$-tame kernel is replaced by its dual.
 Our main theorems provide relationships which do not involve the dual. However, these results are
 for groups of exponent 2. The situation may be different for groups whose $p$-part is non-cyclic for some
 odd prime $p$. We obtain fairly close relationships between $\FI$ and $\SI$, but do not
 investigate the comparison with the Fitting ideal of the dual of $\KEc$, leaving this to the future.

 Of related interest is work of Barrett [\Ba] comparing the image of $\Fit_{\Z[G]}(K_{2k}(\Oc_\E))$
in $\Z_p[G]$, for each odd prime $p$ and positive integer $k$, with
an ideal obtained using values of $L$-functions at the positive
integer $k+1$ and a $p$-adic regulator. Assuming the Equivariant
Tamagawa Number Conjecture and the Quillen-Lichtenbaum conjecture,
he shows that the former ideal contains the latter. This result fits
into the framework of Solomon's conjectures [\So] for $L$-functions
at $s=1$.

\heading{II. Comparisons in Towers of Totally Real Fields}
\endheading

From now on, we assume that $\E$ is a totally real field. Let
$\R=\Z[G]$. We will be considering abelian groups $M$ whose
operation is written multiplicatively, and which possess a natural
$G$-action and therefore become $\R$-modules. For $\alpha\in \R$ and
$m\in M$, we will write $m^\alpha$ for the action of $\alpha$ on
$m$. If $H$ is a subgroup of $G$, then $\Z[H]$ is a subring of $\R$,
and we denote the augmentation ideal in $\Z[H]$ by $I_H$. Set
$\Gb=G/H$. Note that $\R/(I_H)\R\cong \Z[\Gb]$.

\proclaim{Proposition 2.1 (Stickelberger ideals under change of
extension field)} If $E$ is any intermediate field between the
totally real fields $F$ and $\E$, and $\overline{G}=\Gal(E/F)$ then
$\Stick_{E/F}^S(-1)$ equals the image of $\SI$ under the natural
projection map $\pi$ from $\R=\Z[G]$ to
$\overline{\R}=\Z[\overline{G}]$.
\endproclaim

\demo{Proof}
It follows from [\Co, Lemma 2.3] that $\Ann_{\Z[G]}(W_2(\E))$
(respectively $\Ann_{\Z[\overline{G}]}(W_2(E))$) is generated by
all elements of the form $\sigma_\q-\N\q^2$
(respectively $\overline{\sigma}_\q-\N\q^2=\pi(\sigma_\q-\N\q^2)$)
for prime ideals $\q$ not dividing
$w_2(\E)$ and the discriminant of $\E$.  Hence
$\SI=\Ann_{\Z[G]}(W_2(\E))\theta_{\E/F}^S(-1)$
(respectively $\Stick_{E/F}^S(-1) =\Ann_{\Z[\overline{G}]}(W_2(E))\theta_{E/F}^S(-1))$
is generated by
all elements of the form $(\sigma_\q-\N\q^2)\theta_{\E/F}^S(-1)$
(respectively $(\overline{\sigma}_\q-\N\q^2)\theta_{E/F}^S(-1)$).
Extend $\pi$ by $\Q$-linearity;
the inflation property of Artin $L$-functions implies that
$\pi(\theta_{\E/F}^S(-1)) = \theta_{E/F}^S(-1)$.
Hence $\Stick_{E/F}^S$ is generated by the elements
$\pi(\sigma_\q-\N\q^2)\pi(\theta_{\E/F}^S(-1))=\pi\bigl((\sigma_\q-\N\q^2)(\theta_{\E/F}^S(-1))\bigr)$,
and the result follows, since $\pi$ is surjective.
\enddemo

\proclaim{Remark 2.2} {\rm The proof of Proposition 2.1 easily generalizes to higher Stickelberger ideals
for any non-positive integer $-k$. One replaces $\theta_{\E/F}^S(-1)$ by
$\theta_{\E/F}^S(-k)$ and $W_2(\E)$ by $W_{k+1}(\E)$, defined as those roots of unity
fixed by the $k+1$st powers of all automorphisms of $\overline{\Q}$ over $\E$. The Deligne-Ribet theorem
is valid here as well to guarantee that this yields an ideal $\Stick_{\E/F}^S(-k)$ of $\Z[G]$.
This ideal is expected to be related to $\Fit_{\E/F}^S(k)=\Fit_{\Z[G]}(K_{2k}(\OEc))$ (see [\Sn]),
although a closer relationship is
expected upon replacing $K_{2k}(\OEc)$ by an appropriate \'etale cohomology group. In our case of $k=1$,
the two are known to be the same.}
\endproclaim

\proclaim{Proposition 2.3}
Suppose that $\E/E$ is an abelian Galois extension of totally real number fields
with group $H$. Suppose also that $S$ contains all of the infinite primes of
$E$ and the primes which ramify in $\E/E$. Let
$I_H$ denote the augmentation ideal of
$\Z[H]$. Then the transfer map from $\KEc$ to $\KE$ induces an isomorphism
of $\Z[\Gb]$-modules
$$\KEc/\KEc^{I_H}\cong \KE.$$
\endproclaim
\demo{Proof}
Under our assumptions on $S$, it follows more generally from Kahn's
theorem of [\Ka, 5.1] as observed in [\KoMo]
that the transfer map from $\KEc$ to $\KE$ has kernel $\KEc^{I_H}$ and
a cokernel which is an elementary abelian 2-group of rank equal to the
number of real primes of $E$ which ramify in $\E$. Here that rank is zero, and
the result is established.
\enddemo

\proclaim{Proposition 2.4 (Fitting ideals of $S$-tame kernels under change of extension field)}
Suppose that $\E/F$ is an abelian Galois extension of totally real number fields
with group $G$, and $E$ is an intermediate field with $\Gal(E/F)$ denoted by
$\Gb$. Suppose also that $S$ contains all of the infinite primes of
$F$ and the primes which ramify in $\E$.
Then $\FIb$ equals the image of $\FI$ under the natural projection map
$\pi$ from $\R=\Z[G]$ to $\Rb=\Z[\overline{G}]$.
\endproclaim

\demo{Proof} From Proposition 2.3, we get
$$\FIb=\Fit_\Rb(\KE)=\Fit_\Rb(\KEc/\KEc^{I_H}).$$
Identifying $\Rb$ with $\R/(I_H)\R$, this becomes
$$\Fit_{\R/(I_H)\R}(\KEc/\KEc^{I_H})=
\Fit_{\R/(I_H)\R}(\KEc/\KEc^{\R(I_H)}).$$
Now a standard property of Fitting ideals implies that this ideal equals
the image of $\FI=\Fit_{\R}(\KEc)$ in $\R/(I_H)\R$, and under our identifications,
this corresponds to $\pi(\FI)$ in $\Rb$, as desired.
\enddemo

\proclaim{Theorem 2.5 (First Comparison Theorem)}
Suppose that $\E/F$ is an abelian Galois extension of totally real number fields
with group $G$, and $E$ is a relative quadratic extension of $F$ in $\E$.
If the Birch-Tate conjecture holds for $E$ and $F$, then
$\SI$ and $\FI$ in $\Z[G]$ project to the same image in $\Z[\Gal(E/F)]$.
Equality of the images in $\Z[1/2][G]$ holds unconditionally.
\endproclaim

\demo{Proof} Let $\Gb=\Gal(E/F)$ and apply Propositions 2.1 and 2.3.
The images of $\SI$ and $\FI$ in $\Z[\Gb]$ are
$\SIb$ and $\FIb$, respectively. Proposition 1.3 completes the proof.
\enddemo

\proclaim{Proposition 2.6 (Stickelberger ideals under change of base field)}
 Suppose that $\F$ is an intermediate field between $F$ and $\E$. Then
$$\Stick_{\E/\F}^S(-1) \subset \Stick_{\E/F}^S(-1).$$
\endproclaim
\demo{Proof}
This is a special case of the Corollary to the Main Theorem of [\Sa].
\enddemo

\proclaim{Proposition 2.7 (Fitting ideals under change of base field)}
 Suppose that $\F$ is an intermediate field between the totally real fields
 $F$ and $\E$. Then
$$\Fit_{\E/\F}^S(1)\subset \FI.$$
\endproclaim
\demo{Proof}
This is a direct application of a property of Fitting ideals for subrings which
follows immediately from the definition.
\enddemo

\proclaim{Theorem 2.8 (Second Comparison Theorem)}
Suppose that $\F$ is an intermediate field between $F$ and $\E$
such that $\E/\F$ is of degree 2. If the Birch-Tate conjecture holds for $\E$ and $\F$,
then $\FI$ and $\SI$ both contain $\Fit_{\E/\F}^S(1)=\Stick_{\E/\F}^S(-1)$. Without the
assumption of the Birch-Tate conjecture, the extensions of $\FI$ and $\SI$ to $\Z[1/2][G]$
both contain $\Stick_{\E/\F}^S(-1)$.
\endproclaim
\demo{Proof}
Assuming the Birch-Tate conjecture,
the equality $\Fit_{\E/\F}^S(1)=\Stick_{\E/\F}^S(-1)$ holds by Proposition 1.3.
Apply Propositions 2.6 and 2.7 to obtain the result.
Without assuming the Birch-Tate conjecture, the proof goes through
after extending ideals to the group rings obtained by adjoining $1/2$.
\enddemo

\heading{III. Comparisons in the Maximal Order $\S$ of $\Q[G]$ when $G^2=1$}
\endheading

From now on we assume that the abelian group $G$ has exponent 2, and
order $2^m$. Then $\E$ is a composite of relative quadratic
extensions of $F$, and $\E/F$ is what we call a multi-quadratic
extension. In this case, a non-trivial element $\chi$ in the
character group $\Ghat$ of $G$ must have order 2, and a kernel
$\ker(\chi)$ of index 2 in $G$. We denote the fixed field of
$\ker(\chi)$ by $E_\chi$. It is a relative quadratic extension of
$F$. For each non-trivial $\chi$, we fix an element $\tau_\chi\in G$
which does not lie in $\ker(\chi)$. Then $\tau_\chi$ restricts to
the non-trivial automorphism of $E_\chi/F$. Denote the trivial
character of $G$ by $\chi_0$.

Each $\chi$
is associated with an idempotent
$e_\chi=\frac{1}{2^m}\sum_{\sigma \in G}\chi(\sigma)\sigma^{-1} \in \Q[G]$.
The maximal order of
$\R=\Z[G]$ in $\Q[G]$ is $\S=\oplus_{\chi} \Z e_\chi$.
Let $I_\chi$ denote the kernel
of the natural map from $\R$ to $\R e_\chi=\Z e_\chi$. It is easy to see that
$I_\chi$ is generated by the elements $\sigma-\chi(\sigma)$ as
$\sigma$ ranges over $G$.
Note that as $\R$-algebras, we have
$$\S \cong \op_{\chi \in \Ghat} \R e_\chi
\cong \op_\chi \R/I_\chi \R.$$

 To simplify notation, we now set
 $k_2^S(E)=|K_2(\Oc_{E}^S)|$.  For any $\R$-module $M$ and $\alpha \in \R$, let $M_\alpha$ denote the submodule
 of elements annihilated by $\alpha$.
For the relative quadratic extensions $E_\chi/F$, we also set
$k_2^S(E_\chi)^-=|K_2(\Oc_{E_\chi}^S)_{1+\tau_\chi}|$ and
$w_2(E_\chi)^-=|W_2(E_\chi)_{1+\tau_\chi}|$.

\proclaim{Proposition 3.1 (The Stickelberger element in a multi-quadratic extension)}
Assuming that the Birch-Tate conjecture holds for each $E_\chi$ and for $F$, we have
$$(-1)^{|S|}\theta_{\E/F}^S(-1)=\dfrac{\kf}{w_2(F)}e_{\chi_0} +
\sum_{\chi\neq \chi_0}
(-1)^{|S_{E_\chi}|}\dfrac{w_2(F)}{w_2(E_\chi)}\dfrac{\kechi}{\kf}e_\chi.$$
Equality up to factors of 2 in each component holds unconditionally.
\endproclaim
\demo{Proof} This follows from standard properties of $L$-functions
and an application of Proposition 1.2. See [\SaSi, Proposition 5.1]
for more details.
\enddemo

\proclaim{Lemma 3.2} Suppose that $M$ is a finite $\R$-module. Then
$\Fit_{\R}(M)\S = \bigoplus_{\chi \in \Ghat}|M/I_\chi M|\Z e_\chi$.
\endproclaim
\demo{Proof} Note that $\R e_\chi =\Z e_\chi \cong \Z$.
Then using standard properties of Fitting ideals and tensor products, we have
$$
\multline
\Fit_\R(M)\S=\Fit_\S(M\otimes_\R \S)=\Fit_\S(M\otimes_\R (\op_\chi \R/I_\chi))
\\ =\Fit_{\oplus_\chi \R e_\chi}(\op_\chi M/I_\chi M)
=\op_\chi \Fit_{\Z e_\chi}(M/I_\chi M)
\\=\op_\chi |M/I_\chi M|\Z e_\chi.
\endmultline
$$
\enddemo

\proclaim{Lemma 3.3} Suppose that $E$ is an extension of $F$
contained in $\E$ such that $\E/E$ is a quadratic extension. Let
$\tau$ be a generator of $\Gal(\E/E)$. Then as $\Z[G]$-modules,
$\WEc/\WEc^{1-\tau} \cong \WE$ and $\WEc/\WEc^{1+\tau} \cong
\WEc_{1+\tau}$. Also $|W_2(\E)_{1+\tau}| \equiv 2 \pmod{4}$.
\endproclaim
\demo{Proof} From the exact sequence
$$1\rightarrow \WE \rightarrow \WEc \overset{1-\tau}\to{\rightarrow} \WEc^{1-\tau} \rightarrow 1,$$
we see that the first pair of modules to be shown isomorphic indeed have the same order.
Since $\WEc$ is a cyclic group, the isomorphism in question is now clearly induced by raising to the
power $|\WEc^{1-\tau}|$. The proof of the second isomorphism is similar.

For the last statement, clearly $\sqrt{-1}\in W_2(\Q)$ and hence
$\tau$ acts trivially on $\sqrt{-1}$. Thus
$\sqrt{-1}^{1+\tau}=\sqrt{-1}^2=-1\neq 1$. This shows that
$\sqrt{-1}\notin W_2(\E)_{1+\tau}$, so $4 \nmid |W_2(\E)_{1+\tau}|$.
However, $|W_2(\E)_{1+\tau}|$ must be even, as $-1 \in
W_2(\E)_{1+\tau}$.
\enddemo

\proclaim{Corollary 3.4} $\WEc/\WEc^{I_\chi}$ is isomorphic as a $\Z[G]$-module to
$\WF$ when $\chi=\chi_0$ and to $W_2(E_\chi)_{1+\tau_\chi}$ otherwise.
\endproclaim
\demo{Proof} The ideal $I_\chi$ is generated by the elements
$1-\sigma_i$ for $\sigma_i$ running through a minimal set of generators of
$\ker(\chi)$, together with, in the case of $\chi\neq \chi_0$, one element
$1+\sigma_m$ for $\sigma_m=\tau_\chi \notin \ker(\chi)$. Thus the
result follows from successive application of Lemma 3.3 with
$1-\tau=1-\sigma_i$ for $i=1,\dots, m-1$ followed by one more application with $1\pm
\tau =1\pm \sigma_m$.
\enddemo

\proclaim{Proposition 3.5 (Extension of $\Ann_\R(W_2(\E))$ to $\S$)}
$$\Ann_\R{\WEc}\S=w_2(F)\Z e_{\chi_0}\oplus \op_{\chi\neq \chi_0}\wechi^-\Z e_\chi.$$
\endproclaim
\demo{Proof} Since $\WEc$ is a cyclic group, and hence a cyclic $\R$-module,
$$\Ann_\R(\WEc)=\Fit_\R(\WEc).$$ By Lemma 3.2,
$$\Fit_\R(\WEc)\S=\op_{\chi\in \Ghat}|\WEc/\WEc^{I\chi}|\Z e_\chi.$$
Applying Corollary 3.4, we find that
$$\Ann_\R{\WEc}\S=w_2(F)\Z e_{\chi_0}\oplus \op_{\chi\neq \chi_0}|W_2(E_\chi)_{1+\tau_\chi}|\Z e_\chi.$$
\enddemo

For each positive integer $j$, let $F^{(j)}$ denote the $j$th layer
of the cyclotomic $\Z_2$-extension of $F$. It is cyclic of degree
$2^j$ over $F$. Since $G=\Gal(\E/F)$ has exponent 2, $\E$ cannot
contain $F^{(2)}$. For relative quadratic extensions $E/F$ in $\E$,
the case where $E=F^{(1)}$ (if this lies in $\E$), tends to be
rather special. For example, we will make repeated use of the
following result.

\proclaim{Proposition 3.6 (Cohomology of $\WE$ and $K_2(\Oc_{E}^S)$
for a relative quadratic extension)} Suppose that $E/F$ is a
relative quadratic extension of totally real fields, and $\tau$ is
the non-trivial automorphism of $E$ over $F$. Let the finite set $S$
of primes of $F$ contain the infinite and ramified primes as usual.
Let $\iota_{E/F}:K_2(\Oc_{F}^S)\rightarrow K_2(\Oc_{E}^S)_{1-\tau}$
be the natural map induced by the inclusion of rings, with kernel
$\ker(\iota_{E/F})$ and cokernel $\coker(\iota_{E/F})$. Then
$$
\multline
 |\WE_{1-\tau}/\WE^{1+\tau}|=|\WE_{1+\tau}/\WE^{1-\tau}|
\\
=|\KE_{1-\tau}/\KE^{1+\tau}|=|\KE_{1+\tau}/\KE^{1-\tau}|
\\
=|\ker(\iota_{E/F})|=|\coker(\iota_{E/F})|.
\endmultline
$$ If $E =
F^{(1)}$, then the common value of these quantities is 1; otherwise
it is 2. In either case, $\WE_{1-\tau}=\WF$ and
$|\KE_{1-\tau}|=|\KF|$.
\endproclaim
\demo{Proof} See [\San, Propositions 3.4, 2.2 and 2.3].
\enddemo

\proclaim{Lemma 3.7}
 Assuming the Birch-Tate conjecture for $E_\chi$ and $F$, we have
$${w_2(F)}e_{\chi_0}\theta_{\E/F}^S(-1)=\pm \kf e_{\chi_0}$$
and for $\chi \neq \chi_0$,
$$\wechi^-  e_\chi \theta_{\E/F}^S(-1)=\pm \kechi^- e_\chi.$$
Equality of odd parts holds unconditionally.
\endproclaim
\demo{Proof}(Adapted from [\San, Lemma 7.2])
Using Proposition 3.1, we see that
$$
\pm
{w_2(F)}e_{\chi_0}\theta_{\E/F}^S(-1)=w_2(F)\dfrac{\kf}{w_2(F)}e_{\chi_0}={\kf}e_{\chi_0}.
$$
Similarly, and using Lemma 3.3 as well, we have
$$\multline
\pm \wechi^- e_\chi\theta_{\E/F}^S(-1)
\\=|W_2(E_\chi)_{1+\tau_\chi}| \dfrac{w_2(F)}{w_2(E_\chi)}\dfrac{|\KEchi|}{|\KF|}e_\chi
=\dfrac{|W_2(E_\chi)_{1-\tau_\chi}|}{|W_2(E_\chi)^{1+\tau_\chi}|}\dfrac{|\KEchi|}{|\KF|}e_\chi.
\endmultline
$$
At this point, we apply Proposition 3.6 to $E_\chi/F$ and obtain
that $|\KF|=|\KEchi_{1-\tau_\chi}|$, while the cohomology groups of
$\KEchi$ and $W_2(E\chi)$ over $\Gal(E_\chi/F)$ have the same order.
The last expression now becomes
$$\dfrac{|\KEchi_{1-\tau_\chi}|}{|\KEchi^{1+\tau_\chi}|}\dfrac{|\KEchi|}{|\KEchi_{1-\tau_\chi}|}e_\chi
=\dfrac{|\KEchi|}{|\KEchi^{1+\tau_\chi}|}e_\chi=|{K_2(\Oc_{E_\chi}^S)}_{1+{\tau_\chi}}| e_\chi$$
All equalities are valid up to factors of two without the assumption of the Birch-Tate conjecture.
\enddemo

\proclaim{Proposition 3.8 (Extension of the Stickelberger ideal to
$\S$)} Assuming the Birch-Tate conjecture for each $E_\chi$ and for $F$,
$\SI \S $ is generated as an $\S$-ideal by
$$ \kf  e_{\chi_0}+\sum_{\chi \neq \chi_0} \kechi^- e_\chi.$$
Without assuming the Birch-Tate conjecture, this element is a generator for
$\SI \S[1/2]$ in $\S[1/2]=\R[1/2]$.
\endproclaim
\demo{Proof}
Using Proposition 3.5,
$$
\multline
\SI \S = \theta_{\E/F}^S(-1) \Ann_\R(\WEc) \S
\\=\theta_{\E/F}^S(-1) \bigl(w_2(F) e_{\chi_0}+ \sum_{\chi\neq \chi_0} \wechi^- e_\chi \bigr) \S
\endmultline
$$
Since $\S=\oplus{\chi} \Z e_\chi$, it suffices to consider each
component of $\SI \S$ individually. Proposition 3.7 completes the
proof under the assumption of the conjecture. For the unconditional
part, note that the same proof applies in $\S[1/2]=\sum_\chi \Z[1/2]
e_\chi$.
\enddemo

\proclaim{Proposition 3.9 (Extension of the Fitting ideal to $\S$)}
 $$\FI \S = \Z \kf  e_{\chi_0} \oplus \op_{\chi\neq \chi_0} \Z \kechi^- e_\chi$$
\endproclaim
\demo{Proof} Lemma 3.2 gives
 $$\FI \S = \op_{\chi \in \Ghat}|\KEc /\KEc^{I_\chi}| \Z e_\chi.$$
Since $I_{\chi_0}=I_G$, while $I_\chi$ is generated by $I_{\ker(\chi)}$ and $1+\tau_\chi$ for $\chi\neq \chi_0$,
Proposition 2.3 allows us to deduce that
$$\KEc/\KEc^{I_{\chi_0}}\cong \KF$$ and $$\KEc/\KEc^{I_\chi}\cong K_2(\Oc_{E_\chi}^S)/K_2(\Oc_{E_\chi}^S)^{1+\tau_\chi}.$$
Then consideration of the surjective mapping from $K_2(\Oc_{E_\chi}^S)$ to its image under $1+\tau_\chi$
shows that $|K_2(\Oc_{E_\chi}^S)|/|K_2(\Oc_{E_\chi}^S)^{1+\tau_\chi}|=|K_2(\Oc_{E_\chi})^S_{1+\tau_\chi}|$.
This establishes the result.
\enddemo

\proclaim{Theorem 3.10} {\rm (Third Comparison Theorem)} $\FI$ and
$\SI$ extend to the same ideal in $\Z[1/2][G]$. Assuming the
Birch-Tate conjecture for $F$ and for each $E_\chi$, they extend to
the same ideal in $\S$.
\endproclaim
\demo{Proof} The second statement follows directly from Proposition
3.8 and Proposition 3.9. The first follows from these two by taking
images in $\S[1/2]=\R[1/2]=\Z[1/2][G]$.
\enddemo

\proclaim{Lemma 3.11}
Let $p^t$ be a positive power of a prime, and $E$ be a totally real number field.  
Denote a primitive $p^t$th root of unity by $\omega_{p^t}$.
For $p$ odd, we have  $p^t|w_2(E) \iff  \omega_{p^t}+\omega_{p^t}^{-1} \in E$.
For $p=2$, we have $2^{t}|w_2(E) \iff \omega_{2^{t-1}}+\omega_{2^{t-1}}^{-1} \in E$.
\endproclaim 
\demo{Proof}
If $p$ is odd and $p^t|w_2(E)$, then by the definition of $w_2(E)$, 
$\Gal(E(\omega_{p^t})/E)$ has exponent 2. However, this group
contains complex conjugation and injects into the cyclic group $(\Z/p^t\Z)^\times$,
so is generated by complex conjugation. Thus $E$ contains the trace $\omega_{p^t}+\omega_{p^t}^{-1}$.
Conversely, if $E$ contains $\omega_{p^t}+\omega_{p^t}^{-1}$, then $E(\omega_{p^t})/E$ has degree 2,
and $p^t|w_2(E)$.  The proof for $p=2$ requires just a slight modification (see [\SaSi, Lemma 7.2]).
 \enddemo

\proclaim{Corollary 3.12} Let $p^t$ be a positive power of a prime.
\item{1.} If $p^t|w_2(\E)$, then $p^t|w_2(E_{\chi})$ for some non-principal $\chi$.
\item{2.} If $p^t|w_2(\E)$ and $p^{t-1}$ exactly divides $w_2(F)$, then $p|w_2(E_\chi)$ for a unique non-principal $\chi$.
\endproclaim
\demo{Proof}
We assume $p$ is odd. For the proof when $p=2$, simply
replace $t$ by $t-1$ as appropriate, according to Lemma 3.11.
\item{1.} Assume that $p^t|w_2(\E)$. By Lemma 3.11, 
$\omega_{p^t}+\omega_{p^t}^{-1}\in \E$. Thus 
$F(\omega_{p^t}+\omega_{p^t}^{-1})\subset \E$, which is
multi-quadratic over $F$. 
Hence $\Gal(F(\omega_{p^t}+\omega_{p^t}^{-1})/F)$
has exponent 2. However this group also injects into the cyclic group
$(\Z/p^t\Z)^\times /{\pm 1}$, so must have order 1 or 2. So   
$F(\omega_{p^t}+\omega_{p^t}^{-1})$ is contained in a quadratic extension of
$F$ inside $\E$, and hence lies in one of the $E_{\chi}$. By Lemma 3.11 again,
$p^t|w_2(E_{\chi})$.
\item{2.} Let $p^{t-1}$ exactly divide $w_2(F)$. Then $p^{t}|w_2(\E)$ and so by part (1),
$p^{t}|w_2(E_\chi)$ for some $\chi$. By Lemma 3.11, we now have that
$F$ does not contain $\omega_{p^{t}}+\omega_{p^{t}}^{-1}\in E_\chi$.
Since $E_\chi/F$ is relative quadratic, we must have
$E_\chi=F(\omega_{p^{t}}+\omega_{p^{t}}^{-1})$. However, this specifies
$E_\chi$, and hence $\chi$, uniquely. Lemma 3.11 then implies that $p^t|w_2(E_\chi)$
for this $\chi$ and no other. The result follows. 
\enddemo

\proclaim{Lemma 3.13}
Assume the Birch-Tate conjecture holds for $F$, each $E\chi$, and $\E$. Then
$$\frac{|\KEc|/\kf}{w_2(\E)/w_2(F)}=\prod_{\chi\neq \chi_0}\frac{|\KEchi|/\kf}{w_2(E_\chi)/w_2(F)}.$$
Without the assumption of the Birch-Tate conjecture, equality holds for the odd parts.
\endproclaim
\demo{Proof}
Standard properties of Artin $L$-functions give
$$\zeta_\E^S(s)=\zeta_F^S(s)\prod_{\chi\neq \chi_0} L_{\E/F}^S(s,\chi)
=\zeta_F^S(s) \prod_{\chi\neq \chi_0} (\zeta_{E_\chi}^S(s)/\zeta_F^S(s)).$$
Setting $s=-1$ and applying the Birch-Tate conjecture or just its proven odd part
gives the result.
\enddemo

\proclaim{Proposition 3.14 (The index of the higher Stickelberger ideal in a multi-quadratic extension)}
Assume the Birch-Tate conjecture holds for $F$, each $E\chi$, and $\E$. Let
$\delta=\delta_{\E/F}=2$ if $F^{(1)}\subset \E$, and $\delta_{\E/F}=1$ otherwise.
Then
$$(\R : \SI)=|\KEc|\frac{(\SI \S: \SI)}{\delta_{\E/F} 2^{(m-2)2^{m-1}+1}}.$$
Without the assumption of the Birch-Tate conjecture, one has that
$$(\R :\SI)=2^c|\KEc|$$ for some integer $c$.  
\endproclaim
\demo{Proof}
By Proposition 3.8, we have
$$(\S: \SI \S)=\kf \prod_{\chi\neq \chi_0} k_2^S(E_\chi)^-.$$
Applying Proposition 3.6, 
$$\multline
k_2^S(E_\chi)^-=|K_2(\OEchi)_{1+\tau_\chi}|
=|K_2(\OEchi)|/|K_2(\OEchi)^{1+\tau_\chi}|
\\
= \delta_\chi |K_2(\OEchi)|/\kf,
\endmultline$$ 
where $\delta_\chi=1$ if $\E_\chi=F^{(1)}$, and $\delta_\chi=2$ otherwise.
Hence
$$(\S: \SI \S)= (2^{2^m-1}/\delta) \kf \prod_{\chi\neq \chi_0} |K_2(\OEchi)|/\kf.$$
Using Lemma 3.13 then gives us
$$(\S: \SI \S)=(2^{2^m-1}/\delta)|\KE| \frac{\prod_{\chi \neq \chi_0} w_2(E_\chi)/w_2(F)}{w_2(\E)/w_2(F)}.$$
Here, we claim that the term $\frac{\prod_{\chi\neq \chi_0} w_2(E_\chi)/w_2(F)}{w_2(\E)/w_2(F)}$ equals 1.
For, if $p^t$ divides the denominator, then $p^t$ divides the numerator, by Proposition 3.12(1).
On the other hand, if $p^t$ divides the numerator, then it divides only one term in the numerator, by Proposition 3.12(2).
Since $W_2(E_\chi)\subset W_2(\E)$, it is then clear that $p^t$ divides the denominator. We conclude that
$$(\S: \SI \S)=(2^{2^m-1}/\delta)|\KE|.$$ A determinant calculation using the orthogonality relations for characters shows that
$$(\S: \R)=2^{m2^{m-1}}.$$
Combining these clearly yields the result. Without the assumption of the Birch-Tate conjecture, all equalities hold up to powers of 2, and $(\SI \S : \SI)$ is also a power of 2 since $\SI =\SI \R \supset 2^m \SI\S$.   
\enddemo

\heading{IV. Biquadratic Extensions}
\endheading

 Now we assume that $\E/F$ is biquadratic, that is, $m=2$.
 Denote the non-principal characters of $G$ by  $\chi_1$, $\chi_2$, $\chi_3$,
and the corresponding fields by $E_i=E_{\chi_i}$. Also let
$e_i=e_{\chi_i}$ and $\tau_i=\tau_{\chi_i}$, restricting to the
non-trivial automorphism of $E_i/F$. More specifically, we will take
$ \tau_1=\tau_2$ fixing $E_3$ and $\tau_3$ fixing $E_1$.

\proclaim{Proposition 4.1 ($\Ann_\R(\WEc)$ for a biquadratic
extension not containing $F^{(1)}$)} Suppose that $\E/F$ is
biquadratic, and that $F^{(1)}\not\subset \E$. Then
$$
\multline
\Ann_\R(\WEc)=
\\
\Z w_2(F) e_{0} \oplus \op_{1\leq i<j\leq 3} \Z \bigl(w_2(E_i)^- e_i
+ w_2(E_j)^-e_j\bigr)
\\ = (\Ann_\R(\WEc) \S)\cap \R,
\endmultline
$$
of index 2 in $\Ann_\R(\WEc)\S$
\endproclaim
\demo{Proof}
By Proposition 3.5,
$$\bigl(\Ann_\R{\WEc}\bigr)\S=\Z w_2(F) e_{0}\oplus \op_{i} \Z w_2(E_i)^- e_i.$$
This contains $\Ann_\R(\WEc)$, but cannot equal it, as it is not
integral by Lemma 3.3.  We will soon see that $2w_2(E_i)^- e_i \in
\Ann_\R(\WEc)$.

 First let us show that $w_2(F) e_0 \in \Ann_\R(\WEc)$.
Here $e_0=\frac{1}{4}(1+\tau_3)(1+\tau_1)$, as we have specified
that $\tau_3$ fix $E_1$.
 We have seen that $\sqrt{-1}\in W_2(F)$, and thus $4\mid w_2(F)$.
 Hence we may integrally express $w_2(F) e_0= (1+\tau_3)(1+\tau_1)(w_2(F)/4)$.
 Two applications of Proposition 3.6 then yield
 $$
 \multline
 \WEc^{e_0 w_2(F)}=\WEc^{(1+\tau_3)(1+\tau_1)(w_2(F)/4)}
\\ =(W_2(E_1)^2)^{(1+\tau_1)(w_2(F)/4)}
=W_2(E_1)^{(1+\tau_1)(w_2(F)/2)}
\\=(W_2(F)^2)^{w_2(F)/2}=W_2(F)^{w_2(F)}=1
\endmultline
$$

This shows that $w_2(F) e_0 \in \Ann_\R(\WEc)$. It should be clear
that a similar computation using $4e_1=(1+\tau_3)(1-\tau_1)$ will
show that $2w_2(E_1)^- e_1 \in \Ann_\R(\WEc)$, and by the same
token, $2w_2(E_i)^- e_i \in \Ann_\R(\WEc)$.

Finally we check that $w_2(E_1)^- e_1 + w_2(E_2)^- e_2 \in
\Ann_\R(\WEc)$. For this, note that we already know that twice this
element lies in $\Ann_\R(\WEc)$. Hence it suffices to show that this
element is integral and annihilates the 2-part of $\WEc$. We have
taken $\tau_1=\tau_2$ to be the unique non-trivial element fixing
$E_3$, that is, lying in the kernel of $\chi_3$. Since $w_2(E_1)^-
\equiv w_2(E_1)^- \equiv 2 \pmod{4}$ by Lemma 3.3, it follows easily
that the element in question is a $\Z[G]$-multiple of $1-\tau_1$.
Then by Proposition 3.6 applied to $\E/E_3$, we have
$\WEc^{1-\tau_1}=(\WEc_{1+\tau_1})^2$. Lemma 3.3 also implies that
$4 \nmid |(\WEc_{1+\tau_1})|$. Thus the 2-part of
$(\WEc_{1+\tau_1})^2$ is trivial, and this completes the check. Of
course, by symmetry we find that the other elements obtained simply
by permuting the subscripts of $e_1$, $e_2$ and $e_3$ lie in
$\Ann_\R(\WEc)$ as well. Thus the direct sum in the statement of the
proposition lies in $\Ann_R(\WEc)\subset (\Ann_R(\WEc)\S)\cap \R\neq
\Ann_R(\WEc)\S$, and one can easily check that it has index 2 in
$\Ann_\R(\WEc)\S$. At the same time, $\bigl(\Ann_R(\WEc)\S\bigr)
\cap \R$ has index at least 2 in $\Ann_R(\WEc)\S$, and the
conclusion follows.
\enddemo

\proclaim{Corollary 4.2 ($\SI$ for a biquadratic extension not
containing $F^{(1)}$)} Suppose that $\E/F$ is biquadratic, and that
$F^{(1)}\not\subset \E$. Assume that the Birch-Tate conjecture holds
for $F$ and each relative quadratic extension of $F$ in $\E$. Then
$$
\multline \SI
\\=\Z \kf e_{0} \oplus \op_{1\leq i<j\leq 3} \Z
\bigl(k_2^S(E_i)^- e_i + k_2^S(E_j)^- e_j\bigr)
\\ = (\SI \S)\cap (\thE\R),
\endmultline
$$
of index 2 in $\SI \S$.
\endproclaim
\demo{Proof} We begin with the equality in Proposition 4.2 and
multiply by $\thE$ to obtain a formula for $\SI$. Note that $\thE$
is a non-zero divisor in $\Q[G]$ by Proposition 3.1. The result
follows from Lemma 3.7. and the observation that the ambiguity in
sign there affects the generators, but not the $\Z$-module they
generate.
\enddemo

\proclaim{Remark 4.3} {\rm The idempotents in Corollary 4.2 have
denominators equal to 4. However, known results such as [\SaSi, Cor.
6.3, Prop. 6.6] and Proposition 3.6 show that $\kf$ and the
$k_2^S(E_i)^-$ are multiples of 4.}
\endproclaim

  For any prime number $p$, let $\Z_{(p)}\subset \Q$ denote the localization of
$\Z$ at the prime ideal $(p)$. Similarly, if $M$ is a $\Z$-module,
let $M_{(p)}\cong M\otimes_\Z \Z_{(p)}$ be the localization of $M$ at $(p)$.
Note that $\R_{(p)}=\Z_{(p)}[G]\subset \Q[G]$ and the intersection of
these over all primes $p$ is $\R$. Also if $I$ is an ideal of $\R$, then
$I_{(p)}=I \R_{(p)}$ is an ideal of $\R_{(p)}$ and the intersection of these over
all primes $p$ is $I$.

\proclaim{Lemma 4.4} If $$\KEc^{1+\tau_1}\cap\KEc^{1+\tau_3}\neq
\KEc^{(1+\tau_1)(1+\tau_3)}$$ then $\FI$ contains an element which is
congruent to $$k_2^S(E_3)^- e_3\equiv k_2^S(E_3)^-
\frac{1+\tau_1}{2} \pmod{1+\tau_3}.$$
\endproclaim
\demo{Proof} (Compare the proof of [\San, Proposition 9.3]) Since
$4\S \subset \R$, Proposition 3.9 implies that $4 k_2^S(E_3)^- e_3
\in \FI$. Thus $k_2^S(E_3)^- e_3\in \FI_{(p)}$ for each odd prime
$p$. It suffices to show that $\FI_{(2)}$ contains an element
congruent to $k_2^S(E_3)^- e_3 \equiv k_2^S(E_3)^-
\frac{1+\tau_1}{2}$ modulo $(1+\tau_3)$. By the properties of
Fitting ideals, this reduces to considering $M=\KEc/\KEc^{1+\tau_3}$
and showing that $k_2^S(E_3)^- \frac{1+\tau_1}{2}$ is in the Fitting
ideal of $M_{(2)}$ over $\Rb_{(2)}=\R_{(2)}/(1+\tau_3)\cong
\Z_{(2)}[\langle \tau_1 \rangle]$.

  We claim that the cohomology group ${\(M_{(2)}\)}_{1+\tau_1}/M_{(2)}^{1-\tau_1}$
contains an element of order 2 under our hypothesis. We will establish this claim at the
end of the proof. For now, we choose a minimal set of
$\Rb_{(2)}$-generators $\gb_i$ for $M_{(2)}/M_{(2)}^{1-\tau_1}$.
Equivalently, these are a minimal set of generators over
$\Rb_{(2)}/(1-\tau_1)\cong \Z_{(2)}$, that is, generators for this
finite abelian $2$-group. We may assume (see [\San, Lemma 9.2)])
that this group is the direct product of the subgroups of orders
$d_i=2^{c_i}>1$ generated by the $\gb_i$, and that $\gb_1^{d_1/2}\in
\(M_{(2)}\)_{1+\tau_1}/M_{(2)}^{1-\tau_1}$. Then $\prod_i d_i = |
M_{(2)}/M_{(2)}^{1-\tau_1}|$, which is the 2-part of
$|\KEc/\KEc^{(1-\tau_1,1+\tau_3)}|=|K_2(\Oc_{E_3}^S)/K_2(\Oc_{E_3}^S)^{1+\tau_3}|
=|K_2(\Oc_{E_3}^S)_{1+\tau_3}|=k_2^S(E_3)^-$. Hence if we let $D$ be
the diagonal matrix of the $d_i$, then $\det(D)$ is associate to
$k_2^S(E_3)^-$ in $\Z_{(2)}$. We can therefore complete the proof by
showing that $\det(D)\frac{1+\tau_1}{2} \in
\Fit_{\Rb_{(2)}}(M_{(2)})$.

Now $\Rb_{(2)}$ is a local ring and $1-\tau_1$ lies in the maximal
ideal, so by Nakayama's lemma, the (arbitrarily chosen) inverse
images $\gamma_i \in M_{(2)}$ of the $\gb_i$ generate $M_{(2)}$. We
know that $\gamma_i^{d_i}\in M_{(2)}^{1-\tau_1}$; multiplying by
$\gamma_i^{(-d_i/2)(1-\tau_1)}$ shows that
$\gamma_i^{(d_i/2)(1+\tau_1)}\in M_{(2)}^{1-\tau_1}$, for each $i$.
Hence there is a matrix $B$ with entries in $\Rb_2$ such that
$D\bigl(\frac{1+\tau_1}{2}\bigr)-B\bigl({1-\tau_1}\bigr)$ is a
relations matrix for the generators $\gamma_i$ of $M_{(2)}$.
Furthermore, since $\gamma_1^{d_1/2}\in \(M_{(2)}\)_{1+\tau_1}$, we
may choose the first row of $B$ to be zero, so that $\det(B)=0$. The
very definition of the Fitting ideal then gives us that
$$\delta=\det\( D\bigl(\frac{1+\tau_1}{2}\bigr)-B\bigl({1-\tau_1}\bigr)\)\in \Fit_{\Rb_{(2)}}(M_{(2)}).$$
However,
$$
\multline
\delta=\delta\(\frac{1+\tau_1}{2}+\frac{1-\tau_1}{2}\)
 =\delta\frac{1+\tau_1}{2}+\delta\frac{1-\tau_1}{2}
\\ =\det\( D\bigl(\frac{1+\tau_1}{2}\bigr)\)+\det\(-2B \frac{1-\tau_1}{2}\)
\\ =\det(D)\frac{1+\tau_1}{2}+\det(-2B)\frac{1-\tau_1}{2}=\det(D)\frac{1+\tau_1}{2}
\endmultline
$$ as desired.

Finally, we prove the claim. Indeed, we will see that
$|{\(M_{(2)}\)}_{1+\tau_1}/M_{(2)}^{1-\tau_1}|=
|\KEc^{1+\tau_1}\cap\KEc^{1+\tau_3}/ \KEc^{(1+\tau_1)(1+\tau_3)}|$.
For this, note that the cohomology group
${M}_{1+\tau_1}/M^{1-\tau_1}$ over a group of order 2 must have
exponent 2. Thus it is isomorphic to
${\(M_{(2)}\)}_{1+\tau_1}/M_{(2)}^{1-\tau_1}$. We now compute,
using $\KEc/\KEc^{(1-\tau_1)} \cong K_2(\Oc_{E_3}^S)$ from Proposition 2.3.
$$
\multline
|{\(M_{(2)}\)}_{1+\tau_1}/M_{(2)}^{1-\tau_1}|=
|M_{1+\tau_1}/M^{1-\tau_1}|
 = |M/M^{1-\tau_1}|/|M/M_{1+\tau_1}|
\\ =|\KEc/\KEc^{(1-\tau_1,1+\tau_3)}|/|M^{1+\tau_1}|
\\ = |K_2(\Oc_{E_3}^S)/K_2(\Oc_{E_3}^S)^{1+\tau_3}|/|\KEc^{(1+\tau_1,1+\tau_3)}/\KEc^{1+\tau_3}|
\\ = |K_2(\Oc_{E_3}^S)_{1+\tau_3}| / |\KEc^{1+\tau_1}/(\KEc^{1+\tau_1}\cap \KEc^{1+\tau_3})|
\\ = |K_2(\Oc_{E_3}^S)^-|\cdot |\KEc^{1+\tau_3}\cap \KEc^{1+\tau_1}|/|\KEc^{1+\tau_1}|
\endmultline
$$
At this point, let $\epsilon(\E/E_3)=1$ or 2, according to whether
or not $\E=E_3^{(1)}$. Standard properties of the transfer
$\Tr_{\E/E_3}:\KEc\rightarrow K_2(\Oc_{E_3}^S)$ and the natural map
$\iota_{\E/E_3}:K_2(\Oc_{E_3}^S)\rightarrow \KEc$ induced by
inclusion of rings imply that
$\KEc^{1+\tau_1}=\iota_{\E/E_3}(\Tr_{\E/E_3}(\KEc))$. These maps are
also Galois-equivariant. In our case of totally real fields,
$\Tr_{\E/E_3}(\KEc)=K_2(\Oc_{E_3}^S)$ by Proposition 2.3. We
continue our computation with the use of these tools and Proposition
3.6.
$$
\multline
|K_2(\Oc_{E_3}^S)^-|\cdot
|\KEc^{1+\tau_3}\cap \KEc^{1+\tau_1}|/|\KEc^{1+\tau_1}|
\\ =\epsilon(\E/E_3)|K_2(\Oc_{E_3}^S)^-|\cdot |\KEc^{1+\tau_3}\cap \KEc^{1+\tau_1}|/|K_2(\Oc_{E_3})|
\\ =\epsilon(\E/E_3) |\KEc^{1+\tau_3}\cap
\KEc^{1+\tau_1}|/|K_2(\Oc_{E_3})^{1+\tau_3}|
\\
=\frac{|\KEc^{1+\tau_3}\cap\KEc^{1+\tau_1}|}{|\KE^{(1+\tau_1)(1+\tau_3)}|}
\cdot
\frac{\epsilon(\E/E_3)|(\iota_{\E/E_3}(K_2(\Oc_{E_3}^S)))^{1+\tau_3}|}
{|K_2(\Oc_{E_3}^S)^{1+\tau_3}|}
\\ =\frac{|\KEc^{1+\tau_3}\cap\KEc^{1+\tau_1}|}{|\KE^{(1+\tau_1)(1+\tau_3)}|}
\cdot
\frac{\epsilon(\E/E_3)|\iota_{\E/E_3}((K_2(\Oc_{E_3}^S))^{1+\tau_3})|}
{|K_2(\Oc_{E_3}^S)^{1+\tau_3}|}
\\
=\frac{|\KEc^{1+\tau_3}\cap
\KEc^{1+\tau_1}|}{|\KE^{(1+\tau_1)(1+\tau_3)}|}
\cdot
\frac{\epsilon(\E/E_3)}{|\ker(\iota_{\E/E_3}|_{K_2(\Oc_{E_3}^S)^{1+\tau_3}})|}
\endmultline
$$
We now show that the second fraction here is in fact equal to 1. So
consider $|\ker(\iota_{\E/E_3}|_{K_2(\Oc_{E_3}^S)^{1+\tau_3}})|=
|\ker(\iota_{\E/E_3}|_{\iota_{E_3/F}(\KF)})|$. From Proposition 4.1,
we know that $|\ker(\iota_{\E/E_3})|=\epsilon(\E/E_3)\le 2$. When
$\epsilon(\E/E_3)=1$ it is clear that
$|\ker(\iota_{\E/E_3}|_{\iota_{E_3/F}(\KF)})|=1=\epsilon(\E/E_3)$.
When $\epsilon(\E/E_3)=2$, Proposition 7.1 of [\San] shows that a
non-trivial element of $\ker(\iota_{\E/E_3})$ is given by
$\{-1,a\}_{E_3}$, for $a \in E_3$ such that $\E=E_3(\sqrt{a})$.
Since $\E/F$ is biquadratic, we may choose $a \in F$. Then
$\{-1,a\}_{E_3}=\iota_{E_3/F}(\{-1,a\}_F)\in
\ker(\iota_{\E/E_3}|_{\iota_{E_3/F}(\KF)})$, which must then be of
order $2=\epsilon(\E/E_3)$.
\enddemo

\proclaim{Theorem 4.5 (Comparison Theorem for a biquadratic extension not containing
$F^{(1)}$)}
Suppose that $\E/F$ is biquadratic, and that
$\E$ does not contain $F^{(1)}$. Then $\FI\supset 2\FI\S=2\SI \S$.
If  $\KEc^{1+\tau_1}\cap\KEc^{1+\tau_3}\neq
\KEc^{(1+\tau_1)(1+\tau_3)}$ and
 $\KEc^{1+\tau_1}\cap\KEc^{1+\tau_1\tau_3}\neq
\KEc^{(1+\tau_1)(1+\tau_1\tau_3)}$, then either
\item{a.} $\FI=\FI\S=\SI \S \supset \SI$, or
\item{b.} $\FI$ has index 2 in $\FI \S =\SI \S$.
Under the assumption of the Birch-Tate conjecture for
$F$ and the $E_i$, $\FI$ and $\SI$ then have the same index in
$\R$.
\endproclaim
\demo{Proof}
By Proposition 3.9, $\FI\S$ is generated by $|\KF|e_0$
and the $k_2^S(E_i)^- e_i$. We show that twice each of
these elements lies in $\FI$.

Put $\Gb_3=\Gal(E_3/F)$ and consider the projection $\pi_3$ from
$\R=\Z[G]$ to $\Rb=\Z[\Gb_3]$, with kernel generated by $1-\tau_1$.
By Proposition 2.4, $\pi_3(\FI)=\Fit_{E_3/F}^S(1)$. By Proposition
1.3 with $E_3\neq F^{(1)}$, $ \kf \frac{1+\tau_3}{2}$ and
$k_2^S(E_3)^- \frac{1-\tau_3}{2}$ both lie in $\Fit_{E_3/F}^S(1)$.
Thus $\kf \frac{1+\tau_3}{2}+(1-\tau_1)\rho_1$ and $ k_2(E_3)^-
\frac{1-\tau_3}{2}+(1-\tau_1)\rho_2$ lie in $\FI$, for some $\rho_1$
and $\rho_2$ in $\R$. Multiplying by $1+\tau_1$ in $\R$, we deduce
that $2 \kf e_0$ and $2 k_2(E_3)^- e_3$ lie in $\FI$.
 A similar
argument shows that $2 k_2^S(E_2)^- e_2 \in \FI$ and
$2 k_2^S(E_1)^- e_1 \in \FI$. It follows that $2 \FI \S \subset \FI$.

Thus we may consider $\FI/(2\FI\S)$ as an $\FF_2$-subspace of
the 4-dimensional space $(\FI\S)/(2\FI\S)$. If it has dimension 4,
then $\FI=\FI\S$, and this is case (a). If it has dimension 3, then
$\FI$ clearly has index 2 in $\FI \S$. This is case (b).
Under the assumption of the Birch-Tate
conjecture, $\SI$ also has index 2 in
$\FI \S= \SI \S$, by Corollary 4.2 and Theorem 3.10. It follows that $\SI$ and
$\FI$ have the same index in $\R$.

Under our additional assumptions, we now show by contradiction that
$V=\FI/(2\FI\S)$ cannot have dimension less than 3.
For $1\le i \le 3$ let $\Gb_i=\Gal(E_i/F)$ and let $\S_i$ denote the
maximal order in $\Q[\Gb_i]$.
Propositions 2.4 and 1.3 imply that, 
for each $i$ from 1 to 3 inclusive, $V$ projects onto
$\Fit_{E_i/F}^S(1) /
\bigl(2 \Fit_{E_i/F}^S(1) \S_i\bigr)$, of dimension 2, generated by
the images of $\kf e_0$ and $k_2^S(E_i)^- e_i$. So
$V$ has dimension at least 2, hence exactly 2,
and the projection is an isomorphism.
Then $V$ must contain exactly one non-trivial element
whose $e_0$-component is $0$, and the three projections show that
this element must be $v_1= k_2^S(E_1)^- e_1 + k_2^S(E_2)^- e_2 + k_2^S(E_3)^- e_3$.
The subspace $V$ also contains an element with a non-trivial $e_0$-component,
and by adding $v_1$ if necessary, we see that
$V$ contains $v_2=\kf e_0 + k_2^S(E_i)^- e_i$ for some $i$.
For $\sigma_i$ generating $\ker(\chi_i)$ we now take images
modulo $1+\sigma_i$, which amounts to projecting onto the
two-dimensional space spanned by the $e_j$ for $j \neq 0, i$.
Thus the image of $v_2$ modulo $1+\sigma_i$ is 0. The image
of $v_1$ modulo $1+\sigma_i$ has two non-zero components, and thus the image
of $V$ modulo $1+\sigma_i$ does not contain an element
with exactly one non-zero component. This contradicts Lemma 4.4.
\enddemo

\proclaim{Proposition 4.6 ($\Ann_\R(\WEc)$ for a biquadratic
extension containing $F^{(1)}$)} Suppose that $\E/F$ is biquadratic,
and that $E_1=F^{(1)} \subset \E$. Then
$$
\multline \Ann_\R(\WEc)=\Z \bigl(w_2(F) e_{0} + w_2(E_1)^- e_1 +
w_2(E_2)^- e_2\bigr)
\\
\oplus \Z  \bigl(w_2(E_2)^- e_2 + w_2(E_3)^- e_3\bigr)
\\
\oplus \Z 2 w_2(F) e_0 \oplus \Z 2w_2(E_1)^- e_1
\endmultline
$$
of index 4 in $\bigl(\Ann_\R(\WEc)\bigr)\S$
\endproclaim
\demo{Proof} As above, we take $\tau_3$ to fix $E_1$, and
$\tau_1=\tau_2$
 to fix $E_3$. By Proposition 3.5,
$$\Ann_\R(\WEc)\S=\Z w_2(F) e_{0}\oplus \op_{i} \Z w_2(E_i)^- e_i.$$
As in the proof of Proposition 4.1, this contains, but does not
equal $\Ann_\R(\WEc)$, since $w_2(E_i)^- \equiv 2 \pmod{4}$ for each
$i$, by Lemma 3.3. This time, two applications of Proposition
3.6 yield
 $$
 \multline
 \WEc^{e_0 w_2(F)}=\WEc^{(1+\tau_3)(1+\tau_1)(w_2(F)/4)}
\\ =(W_2(E_1)^2)^{(1+\tau_1)(w_2(F)/4)}
=W_2(E_1)^{(1+\tau_1)(w_2(F)/2)}
\\=W_2(F)^{w_2(F)/2}=\{\pm 1\}.
\endmultline
$$
Now we can see that $w_2(F) e_0 \notin \Ann_\R(\WEc)$, but $2 w_2(F)
e_0 \in \Ann_\R(\WEc)$. Similarly, $2 w_2(E_i)^- e_i \in
\Ann_\R(\WEc)$, but $w_2(E_i)^- e_i$ is not integral, so does not
lie in $\Ann_\R(\WEc)$. So far, we know that $2 \Ann_\R(\WEc)\S
\subset \Ann_\R(\WEc)\subset \Ann_\R(\WEc)\S$. The proof that
$w_2(E_2)^- e_2+ w_2(E_3)^-  e_3 \in \Ann_\R(\WEc)$ goes just as in
Proposition 4.2, since we have observed that $F^{(2)}\not \subset
\E$ and thus $\E$ is not the first layer of the cyclotomic
$\Z_2$-extension of $E_1=F^{(1)}$.
 Finally, to see that
$w_2(F) e_{0} + w_2(E_1)^- e_1 + w_2(E_2)^- e_2 \in \Ann_\R(\WEc)$,
note that again we know that twice this element lies in
$\Ann_\R(\WEc)$, and thus is suffices to see that this element
annihilates the 2-part of $\WEc$. Indeed, as $2 w_2(E_i)^-  e_i \in
\Ann_\R(\WEc)$ and $2 w_2(E_i)^-$ is 4 times an odd number for $i\ge
1$, we conclude that $4 e_i$ annihilates the 2-part of $\WEc$ for
$i\ge 1$. After reducing modulo these known annihilators of the
2-part of $\WEc$, it suffices to show that $w_2(F) e_{0} + 2 e_1 + 2
e_2= w_2(F) e_0+ (1-\tau_1)$ annihilates the 2-part of $\WEc$. By
Proposition 3.6 and Lemma 3.3 again,
$\WEc^{1-\tau_1}=W_2(\E)_{1+\tau_1}$, whose 2-part is $\{\pm 1\}$.
Thus if $\omega$ is a generator for the 2-part of $\WEc$, then
$\omega^{1-\tau_1}=-1$. At the same time, we know that $\omega^{e_0
w_2(F)}=-1$ since $\WEc^{e_0 w_2(F)}=\{\pm 1\}$. Combining these
shows that $\omega^{w_2(F) e_0+ (1-\tau_1)}=1$, as desired.

It is easy to check that the ideal in the statement of the
Proposition is contained in $\Ann_\R(\WEc)$ and has index 4 in
$\Ann_\R(\WEc)\S$. If $\Ann_\R(\WEc)$ strictly contains this ideal,
it must be of index 1 or 2 in $\Ann_\R(\WEc)\S$ and lie in $\R$. So
indeed it must equal $\bigl(\Ann_\R(\WEc)\S\bigr) \cap \R$, which is
of index 2 in $\Ann_\R(\WEc)\S$ as in Proposition 4.1. However
$w_2(F) e_0 \in  \bigl(\Ann_\R(\WEc)\S\bigr) \cap \R$, while we have
seen that $w_2(F) e_0 \notin \Ann_\R(\WEc)$. The conclusion follows.
 \enddemo

\proclaim{Corollary 4.7 ($\SI$ for a biquadratic extension
containing $F^{(1)}$)} Suppose that $\E/F$ is biquadratic, and that
$E_1=F^{(1)}$. Assume that the Birch-Tate conjecture holds for $F$
and each relative quadratic extension of $F$ in $\E$. Then
$$
\multline \SI= \Z \bigl(\kf e_{0} + k_2^S(E_1)^- e_1 + k_2^S(E_2)^-
e_2\bigr)
\\
\oplus \Z  \bigl(k_2^S(E_2)^- e_2 + k_2^S(E_3)^- e_3\bigr) \oplus \Z
2 k_2(F) e_0 \oplus \Z 2 k_2^S(E_1)^- e_1
\endmultline
$$
of index 4 in $\SI \S$
\endproclaim
\demo{Proof} We begin with the equality in Proposition 4.6 and
multiply by $\thE$ to obtain a formula for $\SI$. Again, $\thE$ is a
non-zero divisor in $\Q[G]$ by Proposition 3.1. The result follows
from Lemma 3.7 and the observation that the choice of signs there
does not affect the $\Z$-module generated.
\enddemo

\proclaim{Theorem 4.8 (Comparison Theorem for a biquadratic extension
containing $F^{(1)}$)} Suppose that $\E/F$ is biquadratic, and that
$E_1=F^{(1)}$. Then $\FI\supset 2\FI\S=2\SI \S$ and $\FI/(2\FI\S)$ must be
one of three $\FF_2$-subspaces of $\FI\S/(2\FI\S)$ (We cannot say
that all three occur). The bases for these subspaces are:
\item{a.} $\{ \kf e_0 + k_2^S(E_1)^- e_1,\ k_2^S(E_2)^- e_2
,\ k_2^S(E_3)^- e_3\}$
\item{b.} $\{\kf  e_0 +
k_2^S(E_1)^- e_1+ k_2^S(E_2)^- e_2 ,\ k_2^S(E_2)^- e_2+ k_2^S(E_3)^-
e_3\}$
\item{c.} $\{\kf  e_0 +
k_2^S(E_1)^- e_1 ,\ k_2^S(E_2)^- e_2+ k_2^S(E_3)^- e_3\}$

Now assume that the Birch-Tate conjecture holds for $F$ and
each $E_i$.
If case (a) occurs, $\SI$ lies in $\FI$ with index 2. If case (b)
occurs, $\FI=\SI$. If case (c) occurs, $\SI$ and $\FI$ have the same
index in $\R$. If $\KEc^{1+\tau_1}\cap\KEc^{1+\tau_3}\neq
\KEc^{(1+\tau_1)(1+\tau_3)}$, then case (c) does not occur.
\endproclaim
\demo{Proof}
Again by Proposition 3.9, $\FI\S$ is generated by $|\KF|e_0$
and the $k_2^S(E_i)^- e_i$, and we show that twice each of
these elements lies in $\FI$. The proof that
$2 \kf e_0$, $2 k_2(E_3)^- e_3$ and $2 k_2^S(E_2)^- e_2$
lie in $\FI$ goes just as in Theorem 4.5.

For $E_1=F^{(1)}$,
Proposition 1.3 gives $\kf \frac{1+\tau_1}{2}+ k_2^S(E_1)^-
\frac{1-\tau_1}{2}\in \Stick_{E_1/F}^S(-1)$. We obtain $ \kf
\frac{1+\tau_1}{2}+ k_2^S(E_1)^- \frac{1-\tau_1}{2}+(1-\tau_3)\rho_3
\in \FI$, and multiplication by $1+\tau_3$ gives $2 \kf e_0+ 2
k_2^S(E_1)^- e_1 \in \FI$. As we already know that $2 \kf e_0\in
\FI$, we conclude that $2 k_2^S(E_1)^- e_1 \in \FI$

We now consider the images of the $\FF_2$-subspace $\FI/(2\FI\S)$ of
the 4-dimensional space $(\FI\S)/(2\FI\S)$ under projection onto
certain 2-dimensional subspaces. So the kernels of these projections
are 2-dimensional. Let $\S_3$ denote the maximal order in
$\Q[\Gb_3]$. According to Proposition 2.4 and Proposition 1.3,
projecting via $\pi_3$ maps the subspace onto $\Fit_{E_3/F}^S(1) /
\bigl(2 \Fit_{E_3/F}^S(1) \S_3\bigr)$, of dimension 2, generated by
the images of $|\KF|e_0$ and $|K_2(\Oc_{E_3}^S)_{1+\tau_3}|e_3$. So
our subspace in question has dimension at least 2. On the other
hand, projecting via $\pi_1$ yields a one-dimensional image spanned
by the image of $ \kf e_0+ k_2^S(E_1)^-e_1$, according to
Proposition 2.4 and Proposition 1.3, since $E_1=F^{(1)}$. The kernel
has dimension at most 2. So if our subspace has dimension 3, it must
be that it contains the kernel of the projection induced by $\pi_1$,
as well as $\kf e_0+ k_2^S(E_1)^- e_1$. This results in case (a).

 If our subspace is 2-dimensional, it contains just one nontrivial element of
 the kernel of the projection induced by $\pi_1$. If this element
 is not
 $ k_2^S(E_2)^- e_2+ k_2^S(E_3)^- e_3$,
 then the projection of this element via $\pi_2$ or $\pi_3$ will be
 0 and the image of our subspace will be 1-dimensional. We have
 already seen that this is not the case for $\pi_3$, and hence likewise for $\pi_2$.
 So our subspace must contain
 $ k_2^S(E_2)^- e_2+ k_2^S(E_3)^- e_3$,
 which is congruent to 0 modulo $1-\tau_3$,
 and it must also contain an element congruent to $\kf e_0+ k_2^S(E_1)^- e_1$
modulo $1-\tau_3$. This leaves only cases (b) and (c).

The statements concerning $\SI$ are now clear from Corollary 4.7.
Proof of the final claim follows from Lemma 4.4 (and Remark 4.3),
since $\kf e_0\equiv k_2(^S(E_1)^-e_1 \equiv 0$ modulo $1+\tau_3$,
while $k_2^S(E_2)^- e_2 \equiv k_2^S(E_2)^- \frac{1-\tau_1}{2}$.
\enddemo

\proclaim{Proposition 4.9 (The index of the higher Stickelberger ideal for a biquadratic extension)}
Assume that $\E/F$ is biquadratic, and that the Birch-Tate conjecture holds for $\E$, $F$,
and the intermediate fields. Then
$$(\R : \SI) = |\KEc|.$$
\endproclaim
\demo{Proof}
By Proposition 3.14
$$(\R : \SI)=|\KEc|\frac{(\SI \S: \SI)}{2 \delta }.$$
At the same time, Corollaries 4.2 and 4.7 give $(\SI \S : \SI)=2 \delta$.
\enddemo

\proclaim{Remark 4.10 (The index of the higher Stickelberger ideal for a quadratic extension)}
{\rm The result of Proposition 4.9 holds for a relative quadratic extension $\E/F$ as well.
This follows from Proposition 1.3. Still it seems that factors of 2 may intervene in larger 
multi-quadratic extensions.}
\endproclaim

\heading{V. Applications}
\endheading

  For easy reference, we first record some standard facts in a Lemma.

\proclaim{Lemma 5.1} Suppose that $E/F$ is a relative quadratic
extension and that $\alpha$ and $\beta$ lie in $E^\times$. Then
\item{1.}
$E(\sqrt{\alpha})=E(\sqrt{\beta})$ if and only if $\alpha\beta$ is a
square in $E$.
\item{2.} $E(\sqrt{\alpha})/F$ is a Galois extension
if and only if the relative norm of $\alpha$ is a square in $E$.
\item{3.} $E(\sqrt{\alpha})/F$ is a biquadratic extension if and only if
$\alpha$ is not a square in $E$ and the relative norm of $\alpha$ is
a square in $F$.
\endproclaim
\demo{Proof} \item{1.} This follows from Kummer theory or an easy
exercise.
\item{2.} This follows from (1) upon taking $\beta$ to be
the conjugate of $\alpha$ over $F$.
\item{3.} Suppose that the extension is biquadratic. Then
$E(\sqrt{\alpha})=E(\sqrt{a})$ for some $a \in F$. Apply (1) and
take the norm. For the converse, let $c^2$ be the norm of $\alpha$.
The automorphisms sending $\sqrt{\alpha}$ to its conjugates $\pm
c/\sqrt{\alpha}$ both have order two, so cannot lie in a cyclic
group.
\enddemo

\proclaim{Proposition 5.2} Let $F=\Q$ and let $E_1$ be a real
quadratic field of discriminant $d$ for which the prime divisors
$q_j$ of $d$ are not congruent to 1 modulo 4. Let $r$ be a positive,
non-square integer which is a norm from $\Oc_{E_1}$, and whose prime
divisors $p_i$ are also not congruent to 1 modulo 4. Assume further
that $rd$ is not a square. Let $E_3=\Q(\sqrt{r})$, and let
$S$ contain $\{\infty\}\cup \{q_1,\ q_2,\dots \}\cup \{p_1,\ p_2,\dots \}$, but no finite
primes congruent to 1 modulo 4.
Then for $\E=\Q(\sqrt{d},\sqrt{r})$, we have
$\KEc^{1+\tau_1}\cap\KEc^{1+\tau_3}\neq
\KEc^{(1+\tau_1)(1+\tau_3)}$.
\endproclaim
\demo{Proof} Let $\alpha \in \Oc_{E_1}$  have norm $r$. We claim
that the element
 $\{-1,\alpha\}_\E \in \KEc^{1+\tau_1}\cap\KEc^{1+\tau_3}
 =\iota_{\E/E_3}(K_2(\Oc_{E_3}^S))\cap \iota_{\E/E_1}(K_2(\Oc_{E_1}^S))$,
 and $\{-1,\alpha\}_\E \notin \KEc^{(1+\tau_1)(1+\tau_3)}=\iota_{\E/F}(\KF)$.
(These equalities are seen in the proof of Lemma 4.4.)
 First, since $\alpha \in \Oc_{E_1}$ and $N(\alpha)=r$ is an $S$-unit,
 $\alpha$ is also an $S$-unit. Therefore $\{-1,\alpha\}_{E_1}\in K_2(\Oc_{E_1}^S)$
 and clearly $\{-1, \alpha\}_\E=\iota_{\E/E_1}\{-1, \alpha\}_{E_1}\in
 \iota_{\E/E_1}(K_2(\Oc_{E_1}^S))$.

To see that $\{-1,\alpha\}_\E \in \iota_{\E/E_3}(K_2(\Oc_{E_3}^S))$,
 consider the extension $\E(\sqrt{\alpha})/E_3$. The relative
 norm of $\alpha$ in $\E/E_3$ is $r$, which is a square in $E_3$. So this extension
 is biquadratic by Lemma 5.1, and we have $\E(\sqrt{\alpha})=\E(\sqrt{\beta})$
 for some $\beta \in E_3$. Also by Lemma 5.1, $\alpha \beta$ is a
 square in $\E$. It follows that $\{-1,\alpha\}_\E=\{-1,\beta\}_\E
 =\iota_{\E/E_3}(\{-1,\beta\}_{E_3})\in \iota_{\E/E_3}(K_2(\Oc_{E_3}^S))$.
 Note that $\E(\sqrt{\beta})/\Q$ is unramified outside $S\cup
 \{2\}$, and this ensures that $\{-1,\beta\}_{E_3} \in K_2(E_3)$
 actually lies in the $S$-tame kernel
$ K_2(\Oc_{E_3}^S)$ (see [\SaSi, Proposition 6.1]).

 Now we suppose that $\{-1,\alpha\}_\E \in \iota_{\E/\Q}(K_2(\Z^S))$
and derive a contradiction. Being of order 2, this element must be
the image of an element of 2-power order in $K_2(\Z^S)$. Our
choice of $S$ ensures that there are no elements of order 4 in
$K_2(\Z^S)$. For it follows from Tate's computation of $K_2(\Q)$
(see [\Mi, Section 11]) that $K_2(\Z^S)\cong \Z/2\Z \oplus
\bigoplus_{\infty \neq p\in S}(\Z/pZ)^\times$. Thus
$\{-1,\alpha\}_\E$ must be the image of an element of order 2 in
$K_2(\Z^S)$, and such elements are of the form $\{-1,a\}_\Q$ for
some $a\in \Q$, by [\Ta, Theorem 6.1]. We must have
$\{-1,\alpha\}_\E=\{-1,a\}_\E$, or $\{-1,\alpha/a\}_\E=1$. In other
words, $\alpha/a$ is in the Tate kernel. For the totally real field
$\E$, the Tate kernel is generated by $(\E^\times)^2$ and an element
$\pi \in \E$ for which $\sqrt{\pi}$ lies in the cyclotomic
$\Z_2$-extension of $\Q$ (see [\Les,Proposition 2.4]). Hence
$\alpha=a \pi^m \gamma^2$ for some integer $m$ and element
$\gamma\in \E$. So $\E(\sqrt{\alpha})\subset
\E(\sqrt{a},\sqrt{\pi})$, which is an abelian Galois extension of
$\Q$. Consequently, the subfield $E_1(\sqrt{\alpha})$ is also an
abelian Galois extension of $\Q$. On the other hand, the norm of
$\alpha$ is $r$, which is not a square in $\Q$. Furthermore, $r$ is
not a square in $E_1$ for otherwise $E_1=\Q(\sqrt{r})$ and
 $rd$ would be a square in $\Q$
by Lemma 5.1 (1) again. Then by Lemma 5.1 (2), $E_1(\sqrt{\alpha})$
is not a Galois extension of $\Q$, and this is a contradiction.
 \enddemo

\proclaim{Corollary 5.3} Let $r$ be a product of one or more
distinct primes which are congruent to -1 modulo 8, or twice such a
product. Let $S$ contain $\infty$, 2, and the prime divisors of $r$, but no
finite prime congruent to 1 modulo 4.
Then for $F=\Q$ and $\E=\Q(\sqrt{2},\sqrt{r})$, we have
$$\SI \subset \FI,$$
and the index is 1 or 2.
\endproclaim
\demo{Proof} This follows from Proposition 1.2, Proposition 5.2 and Proposition 4.8.
\enddemo

\heading{ACKNOWLEDGEMENTS}
\endheading
We thank Cristian Popescu, David Solomon, and Lloyd Simons for helpful discussions on the subject matter of this paper.
This paper was written during a 2006-2007 sabbatical visit to the Department of Mathematics at the University of California, 
San Diego, whose hospitality we greatly appreciate.

\enddocument
\end